\documentclass[11pt,twoside,a4paper]{article}
\usepackage{amsmath,amssymb,amsthm}
\usepackage{graphicx}
\usepackage[english]{babel}
\usepackage{scrlayer-scrpage}
\usepackage[explicit]{titlesec}
\usepackage{libertine}
\usepackage{tikz}
\usetikzlibrary{calc}
\usepackage{tikzpagenodes}
\usepackage{hyperref}
\usepackage[margin=1.5in]{geometry}

\theoremstyle{plain}
\newtheorem{St}{Theorem}[section]
\newtheorem*{St*}{Theorem}
\newtheorem{Le}[St]{Lemma}
\newtheorem{Gev}[St]{Corollary}
\theoremstyle{definition}

\newtheorem{Def}[St]{Definition}
\newtheorem{Opm}[St]{Remark}

\DeclareMathOperator{\gf}{\mathrm{GF}}
\DeclareMathOperator{\pg}{\mathrm{PG}}
\DeclareMathOperator{\q}{\mathrm{Q}}
\DeclareMathOperator{\w}{\mathrm{W}}
\DeclareMathOperator{\h}{\mathrm{H}}

\newcommand{\cO}{\mathcal{O}}

\begin{document}
	\title{Some non-existence results on $m$-ovoids in classical polar spaces}

\author{Jan De Beule\thanks{Department of Mathematics and Data Science, 
Vrije Universiteit Brussel (VUB),  Pleinlaan 2, B--1050 Brussels, 
Belgium (Email: Jan.De.Beule@vub.be); and Department of 
Mathematics: Logic, Analysis and Discrete Mathematics, Ghent University, Krijgslaan 281 (S8), B-9000 Gent, Belgium.} 
, Jonathan Mannaert\thanks{Department of Mathematics and Data Science, 
Vrije Universiteit Brussel (VUB),  Pleinlaan 2, B--1050 Brussels, 
Belgium (Email: Jonathan.Mannaert@vub.be).}
and Valentino Smaldore\thanks{Dipartimento di Tecnica e Gestione dei Sistemi Industriali,
Universit\`{a} degli Studi di Padova, Stradella S. Nicola 3, 36100 Vicenza, Italy (Email: valentino.smaldore@unipd.it).}}
\maketitle

\begin{abstract}
In this paper we develop non-existence results for $m$-ovoids in the classical polar spaces 
$\q^-(2r+1,q), \w(2r-1,q)$ and $\h(2r,q^2)$ for $r>2$. In \cite{BLP2009} a lower bound on 
$m$ for the existence of $m$-ovoids of $\h(4,q^2)$ is found by using
the connection between $m$-ovoids, two-character sets, and strongly regular graphs. This approach
is generalized in \cite{BKLP} for the polar spaces $\q^-(2r+1,q), \w(2r-1,q)$ and $\h(2r,q^2)$, $r>2$.
In \cite{BDS} an improvement for the particular case $\h(4,q^2)$ is obtained by exploiting the algebraic 
structure of the collinearity graph, and using the characterization of an $m$-ovoid as an intruiging set. 
In this paper, we use an approach based on geometrical and combinatorial arguments, inspired by
the results from \cite{GMP}, to improve the bounds from \cite{BKLP}. 
\end{abstract}

\section{Introduction}

An $m$-ovoid of a polar space is a set $\mathcal{O}$ of points such 
that every generator of the polar space contains exactly $m$ points of $\mathcal{O}$. 
This concept goes back to  Segre, \cite{S}, and it was later introduced for generalized 
quadrangles in \cite{T}, and for polar spaces in \cite{ST}.

In \cite{S}, Segre studied so-called $m$-regular systems of the Hermitian surface $\h(3,q^2)$. By the well 
known fact that the generalized quadrangles $\h(3,q^2)$ and $\q^-(5,q)$ are dually isomorphic,
a $m$-regular system of $\h(3,q^2)$ is an $m$-ovoid of $\q^-(5,q)$. Segre showed that when $q$ is
odd, such an $m$-regular system can only exist if $m=\frac{q+1}{2}$, and called it 
a hemisystem. Moreover, he constructed an example for $q=3$, which gives rise by 
duality to a 2-ovoid of $\q^-(5,3)$. This example was shown to be part of an infinite family for all odd $q$ in 
\cite{CP2005}. In \cite{BH} it was shown that $\q^-(5,q)$ has no $m$-ovoids when $q$ is even.

There are more examples of hemisystems of $\h(3,q^2)$ (and thus
$(\frac{q+1}{2})$-ovoids of $\q^-(5,q)$), see e.g. \cite{BKMX2018,CP2017,CP2015,BGR2013,CP2005,KNS2019, PS2023}.
A second series of constructive results is found in e.g. \cite{K}, where $m'$-ovoids of 
$\h(t(2r+1)-1,q^2)$, $\w(2t(2r+1)-1,q)$ and $\q^{-}(2t(2r+1)-1,q)$ are constructed 
using field reduction, starting from hypothetical $m$-ovoids of the polar space
$\h(2r,q^{2ts})$, see \cite[Corollary 3.3]{K}.

However, similar to $1$-ovoids, $m$-ovoids for $m \geq 1$ in finite classical polar spaces are
very rare objects. Therefore one of the central problems is to determine the values of $m$, 
for which $\mathcal{P}$ has an $m$-ovoid. The main theorem in this paper is the following. 

\begin{St}\label{th:MainImprov}
Let $q > 2$ and $r \geq 3$. Suppose that $\mathcal{O}$ is an $m$-ovoid in one of the following
polar spaces, $Q^-(2r+1,q)$ ($e=2$), $W(2r-1,q)$ ($e=1$) or $H(2r,q)$ ($q$ square, $e=\frac{3}{2})$. 
If (a) $r\geq 4$, or, (b) $e\in\{1,\frac32\}$ and $(r,q,e)\not= (3,3,1)$. Then
$$m\geq\frac{-r(1+\frac{2}{q^{r-e-1}}{+\frac{1}{q^{r-2}}})+\sqrt{r^2(1+\frac{1}{q^{r-e-1}})^2+4(q-2)(r-1)(q^{e+1}\frac{q^{r-2}-1}{q-1}+q^e+1)}}{2(q-1)}.$$
This bound asymptotically converges to $$m\geq \frac{-r+\sqrt{r^2+4(r-1)(q-2)q^{r+e-2}}}{2(q-1)}.$$ 
\end{St}

To prove Theorem~\ref{th:MainImprov}, the techniques found in \cite{GMP}, where a modular 
condition on $m$ for the existence of weighted $m$-ovoids in elliptic quadrics is obtained, have
been inspiring. Theorem~\ref{th:MainImprov} improves asymptotically, of order $\sqrt{r}$, 
the results of \cite[Theorem 13]{BKLP} (Theorem~\ref{th:BKLP}). Improvements
are found for $q \geq 3$ and $r \geq 4$, and subsequent values, and in case of the 
Hermitian polar space $H(2r,q)$ already for $r \geq 3$. A more detailed comparison will be given
in Remark~\ref{Re:improvements}.

%The following result is found in \cite{BKLP}, providing
%a lower bound on $m$ for the existence of $m$-ovoids in the polar spaces $\q^-(2r+1,q), \w(2r-1,q)$ 
%and $\h(2r,q^2)$, $r>2$.

\begin{St}\label{th:BKLP}\cite[Theorem 13]{BKLP}
Consider an $m$-ovoid $\mathcal{O}$ in the polar space $\mathcal{P}_{r,e}'$. 
Then $m \geq b$, with $b$ given in the table below.
\begin{table}[h!]
\begin{center}
\begin{tabular}{|c|c|}
\hline
$\mathcal{P}_{r,e}'$ & $b$ \\
\hline
\rule{0pt}{20pt} $\q^-(2r+1,q)$ & $\displaystyle \frac{-3+\sqrt{9+4q^{r+1}}}{2(q-1)}$ \\[9pt]
\hline
\rule{0pt}{20pt} $\w(2r-1,q)$ & $\displaystyle\frac{-3+\sqrt{9+4q^{r}}}{2(q-1)}$ \\[9pt]
\hline
\rule{0pt}{20pt} $\h(2r,q^2)$ & $\displaystyle\frac{-3+\sqrt{9+4q^{2r+1}}}{2(q^2-1)}$ \\[9pt]
\hline
\end{tabular}
\end{center}
\caption{Lower bounds on $m$.}\label{tab:initial_lower_bounds}
\end{table}
\end{St}

To prove Theorem~\ref{th:BKLP}, the authors first observe that an $m$-ovoid of 
$\mathcal{P}_{r,e}'$ gives rise to a $2$-character set of the ambient projective 
space with relation to the hyperplanes, on its turn inducing a strongly regular graph. The
lower bound on $m$ then follows from the requirement that the parameter 
$\lambda$ of the graph must be non-negative. 

\begin{Opm}\label{rem:one-ovoids}
The lower bound on $m$ in Table~\ref{tab:initial_lower_bounds} is always greater than $1$, except
in the case $\w(2r-1,q)$, $r=2$. Indeed, \cite[Theorem 13]{BKLP} reproduced the non-existence
of $1$-ovoids in the polar spaces $\mathcal{P}'_{r,e}$ (except for $\w(3,q)$) from \cite{Thas1981}.
However, it is well known that $\w(3,q)$ has ovoids if and only if $q$ is even, see e.g. \cite{Tallini}.
\end{Opm}

In \cite{BDS}, the algebraic characterization of an $m$-ovoid as a characteristic vector orthogonal to 
one of the non-trivial eigenspaces of the collinearity graph of the polar space is exploited. 
Also the connection with so-called tight sets is studied in this context. Further combinatorial 
arguments then give the following lower bound on $m$ for $m$-ovoids of $\h(4,q^2)$, $q > 2$.

\begin{St}\cite[Theorem 9.1]{BDS}\label{Th:H(4,q^2)}
Let $\mathcal{O}$ be a non-trivial $m$-ovoid of $H(4,q^2)$. If $q>2$, then
$$m\geq \frac12 \frac{-3q-3+\sqrt{4q^5-4q^4+5q^2-2q+1}}{q^2-q-2}.$$
While for $q=2$, we have $m\geq 2$.
\end{St}

%In \cite{GMP}, a modular condition on $m$ for the existence of weighted $m$-ovoids in elliptic quadrics is obtained. 
%The approach has been inspiring to this paper, enabling us to improve Theorem~\ref{th:BKLP} as follows.

%\begin{St}\label{th:MainImprov}
%Let $q > 2$ and $r \geq 3$. Suppose that $\mathcal{O}$ is an $m$-ovoid in one of the following
%polar spaces, $Q^-(2r+1,q)$ ($e=2$), $W(2r-1,q)$ ($e=1$) or $H(2r,q)$ ($q$ square, $e=\frac{3}{2})$. 
%If (a) $r\geq 4$, or, (b) $e\in\{1,\frac32\}$ and $(r,q,e)\not= (3,3,1)$. Then
%$$m\geq\frac{-r(1+\frac{2}{q^{r-e-1}}{+\frac{1}{q^{r-2}}})+\sqrt{r^2(1+\frac{1}{q^{r-e-1}})^2+4(q-2)(r-1)(q^{e+1}\frac{q^{r-2}-1}{q-%1}+q^e+1)}}{2(q-1)}.$$
%This bound asymptotically converges to $$m\geq \frac{-r+\sqrt{r^2+4(r-1)(q-2)q^{r+e-2}}}{2(q-1)}.$$  
%\end{St}

The main objective of this article is to prove Theorem~\ref{th:MainImprov}. This 
will be achieved by proving Theorem~\ref{Th1:counting} and consequently 
Theorem~\ref{th:MainImprovV1}, solely based on combinatorial and geometrical arguments. 
Theorem~\ref{Th1:counting} also leads to an alternative proof of Theorem~\ref{Th:H(4,q^2)}
and to a generalization to other polar spaces of \cite[Theorem 4.1]{GMP}, see Theorem~\ref{ThSmallImprov}.

%However, using a similar approach, it also can be shown that Theorem~\ref{Th1:counting} can be used to give an easier combinatorial proof to some know results such as Theorem \ref{ThSmallImprov} and Theorem~\ref{Th:H(4,q^2)}. This first theorem essentially reproduce the bounds 
%from Theorem~\ref{th:BKLP} with a tiny improvement, but now
%solely based on combinatorial arguments. The improvement however is to insignificant to consider, but we refer to Remark \ref{Re:th4.3} for more details.
%Finally in Section \ref{sec:newresults} we focus on the main result; i.e. Theorem \ref{th:MainImprov} and a direct consequence i.e. Theorem \ref{th:mainQ7}. Asymptotically it is immediately clear that the results yields improvements, but for the small cases we give some comparisons with older results. These are all denoted in Remark \ref{Re:improvements}.

After the preliminary Section~\ref{sec:prlim} we develop 
some general combinatorial lemmas on weighted $m$-ovoids. In Section~\ref{sec:oldresults} 
we give some proofs of earlier known results, based entirely on combinatorial and 
geometrical arguments. Finally, Section \ref{sec:newresults} is devoted to the proof of 
Theorem~\ref{th:MainImprovV1}, which yields eventually the main result of the paper, 
Theorem~\ref{th:MainImprov}. 

\section{Preliminaries}\label{sec:prlim}

Let $\pg(n,q)$ denote the $n$-dimensional projective space over the finite field $\gf(q)$.
A non-degenerate sesquilinear or non-singular quadratic form on the underlying $(n+1)$-dimensional 
vector space induces a geometry embedded in $\pg(n,q)$, which will be called a {\em finite classical polar space}, or simply
a {\em polar space}. Its elements are the totally isotropic, respectively, totally singular subspaces
with relation to the sesquilinear, respectively quadratic form. The subspaces of maximal dimension contained 
in a polar space are called its {\em generators}, and the projective dimension of a generator is 
$r-1$, with $r$ the Witt index of the underlying form, in which case the polar
space has {\em rank $r$}.

Based on the classification of non-degenerate sesquilinear, respectively, non singular quadratic forms over
finite fields (e.g. \cite{KL1990}), one distinguishes between 3 families of polar spaces: 
orthogonal polar spaces, symplectic polar spaces, and Hermitian polar spaces. Symplectic polar
spaces live only in projective spaces of odd dimension. Orthogonal polar spaces fall apart
in 3 subfamilies, i.e. hyperbolic quadrics, and ellipitic quadrics (both on odd projective
dimension), and parabolic quadrics (in even projective dimension). We distinguish between
Hermitian polar spaces in odd and even projective dimension as well. Note that Hermitian
polar spaces are only defined over a field of square order. 

Let $Q$ be a non-singular quadratic form, then the {\em associated bilinear form $f$} is defined 
as $f(u,v) = Q(u+v) - Q(u) - Q(v)$. In odd characteristic, this relation can be used
to define a non-singular quadratic form from a given non-degenerate orthogonal bilinear form. In
even characteristic, the bilinear form associated to a non-singular quadratic form will be 
a symplectic form, which is non-degenerate if and only if the vector space dimension is even. 
Hence in even characteristic, orthogonal polar spaces are induced only by quadratic forms. 

Hence a polar space $\mathcal{P}$ is either induced by a sesquilinear form, or when it is a 
quadric in characteristic $2$, it has an associated bilinear form. In both cases, this form 
induces a polarity of the ambient projective space, denoted $\perp$, and if $\pi$ is a 
subspace of $\mathcal{P}_{r,e}$, then $\pi^\perp$ is the tangent space to $\mathcal{P}$ 
at $\pi$. Furthermore, for any subspace $\pi$ of the ambient projective space, 
the structure of $\pi^\perp \cap \mathcal{P}$ is well known and this will be 
used throughout this paper.

\begin{table}
\begin{center}
\begin{tabular}{|c|c|c|c|}
\hline
polar space & notation & projective dimension & $e$ \\
\hline
elliptic quadric & $\q^-(2r+1,q)$ & $2r+1$ & 2 \\
hyperbolic quadric & $\q^+(2r-1,q)$ & $2r-1$ & 0 \\
parabolic quadric &  $\q(2r,q)$ & $2r$ & 1 \\
\hline
symplectic space & $\w(2r-1,q)$ & $2r-1$ & 1 \\
\hline
Hermitian polar space &  $\h(2r,q)$ & $2r$ & $3/2$ \\
Hermitian polar space & $\h(2r-1,q)$ & $2r-1$ & $1/2$ \\
\hline
\end{tabular}
\caption{$\mathcal{P}_{r,e}$ polar space of rank $r \geq 1$.}\label{tab:polarspaces}
\end{center}
\end{table}

The notation $\mathcal{P}_{r,e}$ stands for a polar space of rank $r$ 
and type $e$, listed in Table~\ref{tab:polarspaces}. Indeed this notation 
does not distinguish between a symplectic polar space and a parabolic quadric. 
But it will be sufficient to capture all the necessary combinatorial information of the polar space. 

\begin{Def}
Throughout this paper, $\mathcal{P}_{r,e}'$ stands for one of the polar spaces $\w(2r-1,q)$, $\q^-(2r+1,q)$ or $\h(2r,q)$ ($q$ square), i.e.
$e\in \Big\{1, \frac{3}{2}, 2\Big\}$. 
\end{Def}
\begin{Opm}\label{Rem:dim}
Let $\pg(n,q)$ be the ambient projective space of a polar space $\mathcal{P}_{r,e}'$. Then 
\[
n = 2r+2e-3\,.
\]
This will be used throughout the entire paper.
\end{Opm}

\begin{Def}
We use $\theta_n := \frac{q^{n+1}-1}{q-1}$ to denote the number of points in $\pg(n,q)$.
\end{Def}
\begin{Def}
A set $\mathcal{O}$ of points of a polar space $\mathcal{P}_{r,e}$ is an {\em $m$-ovoid 
of $\mathcal{P}_{r,e}$} if and only if every generator of $\mathcal{P}_{r,e}$ contains 
exactly $m$ points of $\mathcal{O}$. It is well known that $|\mathcal{O}|=m(q^{r+e-1}+1)$.
\end{Def}

From the definition it follows immediately that the complement of an $m$-ovoid in the point
set of $\mathcal{P}_{r,e}$ is an $m'$-ovoid of $\mathcal{P}_{r,e}$, with $m' = \theta_{r-1}-m$.

The following lemma characterizes $m$-ovoids of the polar spaces $\mathcal{P}_{r,e}'$.

\begin{Le}(see e.g. \cite[Lemma 1 and Theorem 6]{BKLP}\label{Le:eqdef})
Suppose that $\mathcal{O}$ is a set of points in $\mathcal{P}_{r,e}'$ with ambient projective space $\pg(n,q)$. 
Then the following properties are equivalent.
\begin{enumerate}
\item The set $\mathcal{O}$ is an $m$-ovoid.
\item For every point $p \in \pg(n,q)$
$$|p^\perp \cap \mathcal{O}|=\left\{ \begin{aligned} (m-1)(q^{r+e-2}+1)+1 & , \qquad p\in \mathcal{O}, \\ m(q^{r+e-2}+1) &,\qquad p \in \pg(n,q) \setminus \mathcal{O}\,. \end{aligned}\right.$$
\end{enumerate}
\end{Le}

\begin{Opm}\label{rem:outside}
The previous lemma holds for all finite classical polar spaces under the restriction that $p\in \mathcal{P}_{r,e}$ instead of $\pg(n,q)$.
\end{Opm}

Let $\cO$ be an $m$-ovoid with characteristic vector $\chi$, and let $\pi$ be any subspace of the 
ambient projective space. The {\em weight of $\pi$} is then defined as $\mu(\pi) = \sum_{P \in \pi} \chi_{P}$, 
i.e. the number of points of $\cO$ contained in $\pi$. Hence the characteristic vector of $\cO$ naturally defines 
a weight function on the set of all subspaces of the ambient projective space. 

By Lemma~\ref{Le:eqdef}, a weight function $\mu$ on the subspaces of $\pg(n,q)$, with range $\{0,1\}$ on the
points, defines an $m$-ovoid of $\mathcal{P}_{r,e}'$ if and only if 

\[
\mu(p^\perp)+q^{r+e-2}\mu(p)=m(q^{r+e-2}+1)\,.
\]

This is used in \cite{GMP} to define a {\em(weighted) $m$-ovoid} of a polar space. In the following definition, the domain
of the weight function $\mu$ will be the set of all elements of the polar space, i.e. the points, lines, \ldots, generators.

\begin{Def}\label{def:weightedovoid}
Consider $\mu: \mathcal{P}_{r,e} \rightarrow \mathbb{N}$ such that for every subspace $\pi$ of the ambient projective
space $\pg(n,q)$ it holds that $\mu(\pi)=\sum_{p\in \pi}\mu(p)$. Then we call $\mu$ a 
weighted $m$-ovoid of $\mathcal{P}_{r,e}$ 
if for every point $p$ it holds that
\begin{equation}\label{eq:basic_def}
\mu(p^\perp)+q^{r+e-2}\mu(p)=m(q^{r+e-2}+1)\,.
\end{equation}
\end{Def}

The weight function of a non-weighted ovoid will have the property 
that its range on the points of the polar space is the set $\{0,1\}$.  
We will restrict generally in this paper to non-weighted ovoids, but we 
will use its associated weight function $\mu$ throughout the paper. 
This means that certain results are also valid for weighted $m$-ovoids. 
 
\section{Weighted $m$-ovoids in $\mathcal{P}_{r,e}'$}\label{sec:weightedres}

The main result of \cite{GMP} is a modular condition on $m$ for $m$-ovoids of the elliptic
quadric $\q^-(2r+1,q)$. Some of the combinatorial arguments used to obtain this result,
can be generalized for $m$-ovoids of $\mathcal{P}_{r,e}'$. Recall that $\pg(n,q)$
is the ambient projective space of $\mathcal{P}_{r,e}'$, and $n = 2r+2e-3$.

\begin{Le}\label{Le1}
Suppose that $\mu$ is a weighted $m$-ovoid in $\mathcal{P}_{r,e}'$, 
then for every $j$-dimensional space $\pi$ in $\pg(n,q)$,
$$\mu(\pi^\perp)+q^{r+e-j-2}\mu(\pi)=m(q^{r+e-j-2}+1)\,.$$
\end{Le}
\begin{proof}
Let $A := \{(p,s^\perp) \mid s \in \pi^\perp, p\in s^\perp\}$. We will compute
$$S=\sum_{(p,s^\perp) \in A}\mu(p),$$
using a double counting argument.

Fix the hyperplane $s^\perp$. Each point $p \in s^\perp$ contributes $\mu(p)$ to 
$\mu(s^\perp)$ (using Definition~\ref{def:weightedovoid}). Hence 
using $\dim \pi^\perp=n-j-1$
\begin{equation*}
\begin{split}
S&=\sum_{s\in \pi^\perp} \mu(s^\perp) = \sum_{s\in \pi^\perp} \left(m(q^{r+e-2}+1)-q^{r+e-2}\mu(s)\right)\\
&=  \theta_{n-j-1} m(q^{r+e-2}+1)-q^{r+e-2}\sum_{s\in \pi^\perp}\mu(s)\\
&=  \theta_{n-j-1} m(q^{r+e-2}+1)-q^{r+e-2}\mu(\pi^\perp)\,.\\
\end{split}
\end{equation*}

Fix the point $p$. If $p\in \pi$, then there are $\theta_{n-j-1}$ hyperplanes through $\langle p, \pi\rangle$.
If $p\not\in \pi$, then there are only $\theta_{n-j-2}$. Hence
$$S=\mu(\pi)\theta_{n-j-1}+\left(\mu(\mathcal{P}_{r,e})-\mu(\pi)\right)\theta_{n-j-2}\,.$$

Combining both expressions for $S$ and using $n=2r+2e-3$ proves the statement.
\end{proof}

\begin{St}\label{Th1Gen}\label{Th1:counting}
Suppose that $\mu$ is a weighted $m$-ovoid in $\mathcal{P}_{r,e}'$ and 
let $\pi$ be an arbitrary $j$-dimensional subspace contained in 
$\mathcal{P}_{r,e}'$, $0\leq j\leq r-1$. If $\mu(\pi^\perp\setminus\pi)\neq 0$, 
then
\begin{multline}\label{eq:arbitrarypi}
m(q^{r+e-j-3}+1)(m(q^{r+e-1}+1)-\mu(\pi))+q^{r+e-2}\sum_{p\in \pi^\perp\setminus\pi}\mu(p)^2 = \\
 m(q^{r+e-2}+1)(m-\mu(\pi))(q^{r+e-j-2}+1)+q^{r+e-j-3}\sum_{p\in \mathcal{P}_{r,e}'\setminus \pi}\mu(p)\mu(\langle p,\pi\rangle)+\sum_{s\not\in \pi^\perp}\mu(s^\perp\cap \pi)\,.
\end{multline}

\end{St}
\begin{proof}
\begin{figure}[h!]
    \centering
    \includegraphics[scale=0.2]{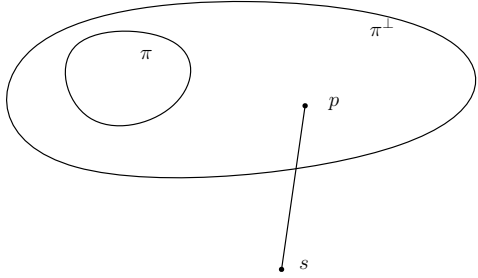}
    \caption{Double counting argument from Theorem \ref{Th1:counting}.}
    \label{fig:Th1:countingl}
\end{figure}
Let $0\leq j\leq r-1$ and fix a $j$-dimensional subspace $\pi$ contained in 
$\mathcal{P}_{r,e}'$. Let $B := \{(p,s)\mid p\in \pi^\perp\setminus\pi, s\not\in\pi^\perp, p \in s^\perp\}$ 
(see Figure \ref{fig:Th1:countingl}). We will compute 
$$S=\sum_{(p,s)\in B} \mu(p)\mu(s).$$
using a double counting argument.

Fix a point $s\not\in \pi^\perp$. Then
\begin{equation*}
\begin{split}
\sum_{p\in s^\perp \cap (\pi^\perp\setminus \pi)} \mu(p) &= \mu(s^\perp \cap (\pi^\perp\setminus \pi))
=\mu(\langle s, \pi \rangle^\perp)-\mu(s^\perp \cap \pi)\\
&=m(q^{r+e-j-3}+1)-q^{r+e-j-3}\mu(\langle s, \pi\rangle)-\mu(s^\perp \cap \pi)\,.
\end{split}
\end{equation*}
This last line follows Lemma \ref{Le1} applied on $\langle s, \pi\rangle$, using $\dim \langle s, \pi \rangle =j+1$.
Hence 
\begin{equation*}
\begin{split}
S&=\sum_{s\not\in \pi^\perp}\mu(s)\left( m(q^{r+e-j-3}+1)-q^{r+e-j-3}\mu(\langle s, \pi\rangle)-\mu(s^\perp \cap \pi)\right)\,.\\
\end{split}
\end{equation*}
Now observe that
\begin{equation*}
\begin{split}
\sum_{s\not\in \pi^\perp}\mu(s)&= \mu(\mathcal{P}_{r,e}\setminus \pi^\perp)
=m(q^{r+e-1}+1)-\left(m(q^{r+e-j-2}+1)-q^{r+e-j-2}\mu(\pi) \right)\\
&=m(q^{r+e-1}-q^{r+e-j-2})+q^{r+e-j-2}\mu(\pi)\,.
\end{split}
\end{equation*}
Using this equation, we find that

\begin{equation*}
\begin{split}
S&= m\left( q^{r+e-j-3}+1\right)\left(m(q^{r+e-1}-q^{r+e-j-2})+\mu(\pi)q^{r+e-j-2}\right) \\&\hspace{4cm}-q^{r+e-j-3}\sum_{s\not\in  \pi^\perp} \mu(s)\mu(\langle s, \pi\rangle)-\sum_{s\not\in \pi^\perp}\mu(s^\perp \cap \pi)\,.
\end{split}
\end{equation*}

Fix a point $p\in \pi^\perp\setminus\pi$ contained in $\mathcal{P}_{r,e}'$. 
Then

$$\sum_{s\in p^\perp, s\not \in \pi^\perp}\mu(s)= \mu(p^\perp)-\mu(\langle p, \pi \rangle^\perp).$$
Hence 
\begin{equation*}
\begin{split}
S&=\sum_{(p,s)\in B} \mu(p)\mu(s)
=\underbrace{\sum_{p\in \pi^\perp\setminus\pi} \mu(p^\perp)\mu(p)}_{:=S_1}-\underbrace{\sum_{p\in \pi^\perp\setminus\pi} \mu(\langle p,\pi \rangle^\perp)\mu(p)}_{:=S_2} \,.
\end{split}
\end{equation*}
By Lemma~\ref{Le1}, applied on $p$, we find 
\begin{equation*}
\begin{split}
S_1&=\sum_{p\in \pi^\perp\setminus\pi} \mu(p^\perp)\mu(p)
= \sum_{p\in  \pi^\perp\setminus\pi} (m(q^{r+e-2}+1)-q^{r+e-2}\mu(p))\mu(p)\\
&=m(q^{r+e-2}+1)\mu( \pi^\perp\setminus\pi)-q^{r+e-2} \sum_{p\in  \pi^\perp\setminus\pi} \mu(p)^2\,,
\end{split}
\end{equation*}
and, applied on $\pi$,
\begin{equation*}
\begin{split}
\mu( \pi^\perp\setminus\pi)  = \mu(\pi^\perp)-\mu(\pi)
&=m(q^{r+e-j-2}+1)-q^{r+e-j-2}\mu(\pi)-\mu(\pi)\\
 & = (m-\mu(\pi))(q^{r+e-j-2}+1)\,.
\end{split}
\end{equation*}
Hence
$$S_1=m(q^{r+e-2}+1)(m-\mu(\pi))(q^{r+e-j-2}+1)-q^{r+e-2} \sum_{p\in  \pi^\perp\setminus\pi} \mu(p)^2.$$

Similarly, once again applying Lemma~\ref{Le1}, we find 
\begin{equation*}
\begin{split}
S_2&=\sum_{p\in \pi^\perp\setminus\pi} \mu(\langle p,\pi \rangle^\perp)\mu(p)\\
&= \sum_{p\in \pi^\perp\setminus\pi} \left(m(q^{r+e-j-3}+1)-q^{r+e-j-3}\mu(\langle p, \pi \rangle)\right)\mu(p)\\
&=m(q^{r+e-j-3}+1)\mu(\pi^\perp\setminus\pi)-q^{r+e-j-3}\sum_{p\in\pi^\perp\setminus\pi}\mu(p)\mu(\langle p, \pi\rangle)\\
&=m(q^{r+e-j-3}+1)(m-\mu(\pi))(q^{r+e-j-2}+1)-q^{r+e-j-3}\sum_{p\in \pi^\perp\setminus\pi}\mu(p)\mu(\langle p, \pi\rangle)
\end{split}
\end{equation*}
The statement of the Lemma now follows from both equations for $S = S_1 - S_2$,

$$m(q^{r+e-j-3}+1)(m(q^{r+e-1}+1)-\mu(\pi))+q^{r+e-2}\sum_{p\in \pi^\perp\setminus\pi}\mu(p)^2$$
$$=m(q^{r+e-2}+1)(m-\mu(\pi))(q^{r+e-j-2}+1)+q^{r+e-j-3}\sum_{p\in \mathcal{P}_{r,e}'\setminus \pi}\mu(p)\mu(\langle p,\pi\rangle)+\sum_{s\not\in \pi^\perp}\mu(s^\perp\cap \pi).$$

\end{proof}

In the particular case where $j=0$, we find the following equation.

\begin{Gev}\label{Th1}\label{Cor:countingpoint}
Suppose that $\mu$ is a weighted $m$-ovoid in $\mathcal{P}_{r,e}'$ and 
let $p_0$ be an arbitrary point in $\mathcal{P}'_{r,e}$ such that $\mu(p_0)<m$. Then
\begin{equation*}
\begin{split}
m(q^{r+e-3}+1)(m(q^{r+e-1}+1)-\mu(p_0))+ q^{r+e-2} \sum_{p\in p_0^\perp\setminus\{p_0\}} \mu(p)^2\\
=m(q^{r+e-2}+1)^2(m-\mu(p_0))+q^{r+e-3}\sum_{p\in \mathcal{P}_{r,e}'\setminus\{p_0\}}\mu(p)\mu(\langle p_0, p\rangle)
\end{split}
\end{equation*}
\end{Gev}
\begin{proof}
If $\mu(p_0) < m$, then it follows immediately from \eqref{eq:basic_def} (Definition~\ref{def:weightedovoid}) that 
$\mu(p_0^\perp\setminus\{p_0\})\neq 0$. Now we can use Theorem~\ref{Th1:counting}.
\end{proof}

\section{A lower bound on $m$ for $m$-ovoids in $\mathcal{P}'_{r,e}$}\label{sec:oldresults}

Theorem \ref{Th1:counting} and Corollary \ref{Cor:countingpoint} will show to be very useful in studying non-weighted $m$-ovoids. 
In this section, the non-weighted ovoid is represented by its weight function $\mu$, which is now a $\{0,1\}$- valued function 
when restricted to the points of the polar space, and even the ambient projective space. Also note that 
an $m$-ovoid $\cO$ of $\mathcal{P}'_{r,e}$ induces a unique weight function $\mu$, and vice versa. 
We call an $m$-ovoid non-trivial if $m \neq 0$ and $m$ does not equal the number of points in a generator.

\begin{Le}\label{Le2}
Let $\cO$ be a non-trivial $m$-ovoid in $\mathcal{P}_{r,e}'$, with weight function $\mu$. 
\begin{itemize}
\item If $p_0\in \mathcal{P}_{r,e}'$, then $\sum_{p\in p_0^\perp\setminus\{p_0\}} \mu(p)^2=(m-\mu(p_0))(q^{r+e-2}+1)$,
\item and if $p_0\in \mathcal{O}$ then $$\sum_{p\in \mathcal{P}_{r,e}'\setminus\{p_0\}}\mu(p)\mu(\langle p_0, p\rangle)\geq 2(m(q^{r+e-1}+1)-1)\,.$$
\end{itemize}
\end{Le}
\begin{proof}
Note that $\sum_{p\in p_0^\perp\setminus\{p_0\}} \mu(p)^2=(m-\mu(p_0))(q^{r+e-2}+1)$ by considering all generators of $\mathcal{P}_{r,e}'$
through $p_0$, and that $\mu(p)^2 = \mu(p)$ since $\mu$ is $\{0,1\}$-valued on the points of $\mathcal{P}_{r,e}'$. 
The first point follows from the fact that $\mu$ is a $\{0,1\}$-valued function. 

Now let $p \in \cO \setminus \{p_0\}$. Then $\mu(\langle p, p_0\rangle)\geq \mu(p)+\mu(p_0)=2$, and this gives now immediately the statement. 
\end{proof}

Now we can substitute the bounds of Lemma~\ref{Le2} in the equation from Corollary~\ref{Cor:countingpoint}. This gives the following theorem.

\begin{St}\label{th:equation}
Let $\mathcal{O}$ be a non-trivial $m$-ovoid of $\mathcal{P}_{r,e}'$, 
with $m \geq 2$. then 
\begin{equation}\label{eq:mainequationofthissection}
(q-1)^2m^2+3(q-1)m-q^{r+e-1}-q-2\geq 0.
\end{equation}
\end{St}
\begin{proof}
%If all points of $\cO$ would have weight $m$, then $\mu' := \frac{1}{m} \mu$ is a $\{0,1\}$-valued 
%function on the points of $\cO$ inducing a $1$-ovoid. So we may assume that not all points have weight
%$m$. 
Since $m \geq 2$, we can find a point $p_0 \in \cO$ with $\mu(p_0) < m$. For this point $p_0$
we can apply Corollary~\ref{Cor:countingpoint} to get
\begin{equation*}
\begin{split}
m(q^{r+e-3}+1)(m(q^{r+e-1}+1)-\mu(p_0))+ q^{r+e-2} \sum_{p\in p_0^\perp\setminus\{p_0\}} \mu(p)^2\\
=m(q^{r+e-2}+1)^2(m-\mu(p_0))+q^{r+e-3}\sum_{p\in \mathcal{P}_{r,e}'\setminus\{p_0\}}\mu(p)\mu(\langle p_0, p\rangle).
\end{split}
\end{equation*}
For the same point $p_0$ we can use the bounds of Lemma~\ref{Le2}, and find
\begin{equation*}
\begin{split}
m(q^{r+e-3}+1)(m(q^{r+e-1}-1)+q^{r+e-2}(m-\mu(p_0))(q^{r+e-2}+1)\\ 
\geq m(q^{r+e-2}+1)^2(m-\mu(p_0))+2q^{r+e-3}(m(q^{r+e-1}+1)-1)\,.
\end{split}
\end{equation*}
This equation simplifies to
\[
q^{r+e-3}(q-1)^2m^2+3q^{r+e-3}(q-1)m+q^{r+e-3}(-q^{r+e-1}-q+2)\geq 0.
\]Dividing by $q^{r+e-3}$ gives the desired inequality.
\end{proof}

The following theorem was also proven in \cite[Theorem 4.1]{GMP} for $\q^-(2r+1,q)$ and can now be generalized for any of the polar spaces $\mathcal{P}_{r,e}'$.
\begin{St}\label{ThSmallImprov}
Consider a non-trivial $m$-ovoid $\mathcal{O}$ in $\mathcal{P}_{r,e}'$, $r \geq 2$.
\begin{enumerate}
\item If $\mathcal{P}_{r,e}'= \q^-(2r+1,q)$ then
$$ m\geq \frac{-3+\sqrt{9+4(q^{r+1}+q-2)}}{2(q-1)}.$$
\item If $\mathcal{P}_{r,e}'= \w(2r-1,q)$ and $r > 2$ then
$$m\geq \frac{-3+\sqrt{9+4(q^{r}+q-2)}}{2(q-1)}.$$
\item If $\mathcal{P}_{r,e}'=\h(2r,q^2)$ then
$$ m\geq \frac{-3+\sqrt{9+4(q^{2r+1}+q^2-2)}}{2(q^2-1)}.$$
\end{enumerate}
\end{St}
\begin{proof}
We may assume that $m > 1$, since $1$-ovoids are excluded for the polar spaces under consideration, 
see Remark~\ref{rem:one-ovoids}. So we can apply Theorem~\ref{th:equation}. We consider 
the three different cases. It is sufficient to fill in the parameter $e$ of the 
particular polar space and to use the correct order of the underlying field in 
Equation~\eqref{eq:mainequationofthissection}. We always obtain a positive lower bound on $m$.

Let $e=2$, then $\mathcal{P}_{r,e}'=\q^-(2r+1,q)$. Then Equation~\eqref{eq:mainequationofthissection} becomes 
$$(q-1)^2m^2+3(q-1)m-q^{r+1}-q+2\geq 2\,.$$
It follows that 
$$m\leq \frac{-3-\sqrt{9+4(q^{r+1}+q-2)}}{2(q-1)} \text{ or } m\geq \frac{-3+\sqrt{9+4(q^{r+1}+q-2)}}{2(q-1)}.$$

Let $e=1$, then $\mathcal{P}_{r,e}' = \w(2r-1,q)$, Then Equation~\eqref{eq:mainequationofthissection} becomes 
$$(q-1)^2m^2+3(q-1)m-q^{r}-q+1\geq 0.$$
Similarly, it follows that 
$$m\leq \frac{-3-\sqrt{9+4(q^{r}+q-2)}}{2(q-1)} \text{ or } m\geq \frac{-3+\sqrt{9+4(q^{r}+q-2)}}{2(q-1)}.$$

Finally let $e=\frac{3}{2}$, then $\mathcal{P}_{r,e}'=\h(2r,q^2)$. Then Equation~\eqref{eq:mainequationofthissection} becomes 
$$(q^2-1)^2m^2+3(q^2-1)m-q^{2r+1}-q^2+2\geq 0.$$
Similarly, it follows that 
$$m\leq \frac{-3-\sqrt{9+4(q^{2r+1}+q^2-2)}}{2(q^2-1)} \text{ or } m\geq \frac{-3+\sqrt{9+4(q^{2r+1}+q^2-2)}}{2(q^2-1)}.$$

\end{proof}

\begin{Opm}\label{Re:th4.3}
\begin{itemize}
\item The bounds of Theorem~\ref{ThSmallImprov} slightly improve the bounds of Theorem~\ref{th:BKLP} by adding a term $4(s-2)$ under the square root, with $s$ the order
of the underlying field.
\item For the symplectic case, it is easy to check that $\frac{-3+\sqrt{9+4(q^{r}+q-2)}}{2(q-1)} = 1$ for $r=2$. Hence this bound is obviously 
also correct, but cannot be improved anyway, since $\w(3,q)$ has ovoids if and only if $q$ is even, see also Remark~\ref{rem:one-ovoids}.
\end{itemize}
\end{Opm}

\subsection*{A lower bound on $m$ for $m$-ovoids of $\h(4,q^2)$}

The generators of the polar space $\h(4,q^2)$ are lines. This enables us to show a stronger version of Lemma~\ref{Le2},
which in its turn can be used to show Theorem~\ref{Th:H(4,q^2)} ($q>2$). As explained in the introduction, we thus 
obtain a proof of Theorem~\ref{Th:H(4,q^2)} ($q>2$) which only relies on combinatorial and geometrical arguments, and not
on the fact that an $m$-ovoids meets an $i$-tight set in exactly $mi$ points, which is proven in \cite{BDS} 
for generalized quadrangles of order $(s^2,s^3)$ (and hence also for $\h(4,q^2)$) using the 
underlying algebraic properties of the collinearity graph.

\begin{Le}\label{Le:UpdateH(4,q^2)}
Let $\mathcal{O}$ be a non-trivial $m$-ovoid in $\h(4,q^2)$ with weight function $\mu$. Fix 
a point $p_0\in \h(4,q^2)\cap \mathcal{O}$, then 
\begin{equation}\label{eq:findanoriginallabel}
\sum_{p\in H(4,q^2)\setminus\{p_0\}}\mu(p)\mu(\langle p_0, p\rangle)\geq m(m-1)(q^3+1)+2(mq^3(q^2-1)+q^3)\,.
\end{equation}
\end{Le}
\begin{proof}
Consider the point $p_0$ and consider all lines through this point in $p_0^\perp$. 
These lines are tangent lines to or generators of $\h(4,q^2)$. There are $q^3+1$ generators 
through $p_0$. Each such generator meets $\cO$ in exactly $m$ points. Hence on each generator
through $p_0$ there are exactly $m-1$ points of $\cO \setminus \{p_0\}$. So for a point $p \in \cO \cap p_0^\perp$,
$\mu(\langle p_0,p \rangle) = m$. For any other point $p \in \cO \setminus p_0^\perp$, $\mu(\langle p_0,p \rangle) \geq 2$.
The number of points in $\cO \setminus p_0^\perp$ equals 
$$m(q^5+1)-1 = m(q^5+1)-(m-1)(q^3+1)-1=2(mq^3(q^2-1)+q^3),$$
and now similarly as in the proof of Lemma~\ref{Le2}, we can conclude Inequality~\eqref{eq:findanoriginallabel}.
\end{proof}

Using this improvement of Lemma \ref{Le2} we can show the following theorem.

\begin{St}\cite[Theorem 9.1]{BDS}
Let $\mathcal{O}$ be an $m$-ovoid of $\h(4,q^2)$, for $q>2$, then
$$m\geq \frac{-3q-3+\sqrt{4q^5-4q^4+5q^2-2q+1}}{2(q^2-q-2)}\,.$$
\end{St}
\begin{proof}
Denote as usual by $\mu$ the weight function of $\cO$. Fix a point $p_0\in \mathcal{O}$.
By Corollary~\ref{Cor:countingpoint}, Lemma \ref{Le2} and Lemma~\ref{Le:UpdateH(4,q^2)}
we find the following inequality.
\begin{equation*}
\begin{split}
m(q+1)(m(q^5+1)-1)&+q^3(q^3+1)(m-1)\\ &\geq m(q^3+1)^2(m-1)+q(q^3+1)m(m-1)+2q(mq^3(q^2-1)+q^3).
\end{split}
\end{equation*}
This inequality simplifies to 
$$(q^5-q^4-2q^3)m^2+3q^3(q+1)m-q^6-q^3-2q^4\geq 0,$$
and further to 
$$(q^2-q-2)m^2+3(q+1)m-q^3-1-2q\geq 0,$$
hence
$$
m(m(q-2)+3) \geq \frac{q^3+2q+1}{q+1}\,.
$$
Hence for $q > 2$ and the fact that $m$ has to be positive, the statement of the lemma follows.
%Note that when $q>2$ the factor $q^5-q^4-2q^3>0$, hence we obtain that
%\begin{equation*}
%\begin{split}
%m&\leq \frac{-3q-3-\sqrt{4q^5-4q^4+5q^2-2q+1}}{2(q^2-q-2)}\\
%\text{ or }\\
%m&\geq \frac{-3q-3+\sqrt{4q^5-4q^4+5q^2-2q+1}}{2(q^2-q-2)}\\
%\end{split}
%\end{equation*}
\end{proof}

\section{New non-existence conditions on $m$-ovoids}\label{sec:newresults}
In this particular section, we focus on new non-existence results for $m$-ovoids in the polar spaces $\mathcal{P}'_{r,e}$
essentially using Theorem \ref{Th1Gen}. %But first we need the following results.

\begin{Le}\label{Le:aid1}
Let $\mathcal{O}$ be an $m$-ovoid in $\mathcal{P}'_{r,e}$ with weight function $\mu$. 
If $\pi$ is an $(r-2)$-space contained in $\mathcal{P}'_{r,e}$, then 
$$\sum_{s\not\in \pi^\perp}\mu(s^\perp\cap \pi)= \mu(\pi)q^{r+2e-1}\frac{q^{r-2}-1}{q-1}.$$
\end{Le}
\begin{proof}
%Fix an $(r-2)$-dimensional subspace $\pi$ contained in  $\mathcal{P}'_{r,e}$. 
We double count the number of pairs in the 
set $C:=\{(p,s) \mid s\not\in \pi^\perp, p\in s^\perp\cap \pi\cap\mathcal{O}\}$. Fix a point $s \not \in \pi^\perp$. 
Then clearly the number of points $p \in s^\perp \cap \pi \cap \cO$ equals $\mu(s^\perp \cap \pi)$. 
Hence $|C| = \sum_{s\not\in \pi^\perp}\mu(s^\perp\cap \pi)$. Now fix a point $p \in \pi \cap \cO$. The number of
points $s$ equals the number of points in $p^\perp \setminus \pi^\perp$, which is
$$|p^\perp|-|\pi^\perp|=\frac{q^{2r+2e-3}-1}{q-1}-\frac{q^{r+2e-1}-1}{q-1}=q^{r+2e-1}\frac{q^{r-2}-1}{q-1}.$$
Hence $|C| = \sum_{p \in \pi \cap \cO} q^{r+2e-1}\frac{q^{r-2}-1}{q-1} = \mu(\pi)q^{r+2e-1}\frac{q^{r-2}-1}{q-1}$.
Combining both expressions for $|C|$ gives the statement of the lemma.
\end{proof}

\begin{Le}\label{Le:aid2}
Let $\mathcal{O}$ be an $m$-ovoid in $\mathcal{P}'_{r,e}$ with weight function $\mu$. 
If $\pi$ is an $(r-2)$-space contained in $\mathcal{P}'_{r,e}$, then 
$$\sum_{p\in \mathcal{P}_{r,e}'\setminus \pi}\mu(p)\mu(\langle p,\pi\rangle)\geq m(q^e+1)(m-\mu(\pi))+(1+\mu(\pi))(mq^e(q^{r-1}-1)+\mu(\pi)(q^e+1)).$$
\end{Le}
\begin{proof}
%Fix an $(r-2)$-dimensional subspace of contained in $\mathcal{P}'_{r,e}$. 
Recall that $\pi^\perp \cap \mathcal{P}'_{r,e}$ is a cone with vertex $\pi$ and base $\mathcal{P}'_{1,e}$. All
generators through $\pi$ are contained in $\pi^\perp$, and each generator meets $\cO$ in exactly $m$ points. 
Hence $\mu(\pi^\perp\setminus \pi) = (m-\mu(\pi))(q^e+1)$. If a point $p \in \pi^\perp\setminus \pi$, then $\langle p , \pi \rangle$ is 
a generator and $\mu(\langle p, \pi \rangle) = m$, while for a point $p \not \in \pi^\perp$, 
$\mu(\langle p, \pi \rangle) \geq \mu(p) + \mu(\pi)$. Hence if $p\in \pi^\perp \setminus \{\pi\}$, then 
$$\sum_{p\in \pi^\perp\setminus \pi}\mu(p)\mu(\langle p,\pi\rangle)= m \mu(\pi^\perp\setminus\pi) = m(m-\mu(\pi))(q^e+1)\,,$$
and if $p\not \in \pi^\perp$, then
\begin{equation*}
\begin{split}
\sum_{p\not \in\pi^\perp}\mu(p)\mu(\langle p,\pi\rangle) \geq &  \mu(\mathcal{P}'_{r,e}\setminus\pi^\perp)(\mu(\pi)+\mu(p)) \\
 &  = (m(q^{r+e-1}+1)-(m-\mu(\pi))(q^e+1))(1+\mu(\pi)).
\end{split}
\end{equation*}
Summing up both expressions yields the desired inequality. 
\end{proof}
%Using these two results will yield the following observation.

\begin{Le}\label{Le:aid3}
    Suppose that $\mathcal{O}$ is an $m$-ovoid in $\mathcal{P}'_{r,e}$, then there exist an $(r-2)$-space with at least $\min\{m, r-1\}$ points of $\mathcal{O}$.
\end{Le}
\begin{proof}
Let $\tau$ be a generator, then $\tau$ meets $\cO$ in exactly $m$ points. 
Let $n = \min \{m,r-1\}$, and let $A$ be any subset of $\tau \cap \cO$ of size $n$. 
Then $A$ cannot span a generator, so $\langle A \rangle$ is a subspace contained in $\tau$ of 
dimension at most $r-2$ containing at least $\min \{m,r-1\}$ points of $\cO$. 
\end{proof}
\begin{St}\label{th:MainImprovV1}
Assume that $\mathcal{O}$ is an $m$-ovoid in $\mathcal{P}_{r,e}'$ and that $\pi$ is an arbitrary $(r-2)$-space contained 
in $\mathcal{P}_{r,e}'$  such that $\mu(\pi^\perp\setminus\{\pi\})\neq 0$, then
\begin{equation}\label{EqNew}
\begin{split}
m^2(q^{r}&-q^{r-1}-q^{e}-q) +m\left( \mu(\pi)(q^{r-1}+2q^{e}+q)+q^{r-1}+q^{e}\right)\\
&-\mu(\pi)\left( q^{r+e-1}+q^{r-1}+(1+\mu(\pi))(q^{e}+1)+q^{r+e}\frac{q^{r-2}-1}{q-1}\right)\geq 0.
\end{split}
\end{equation}
\end{St}
\begin{proof}
Consider the $(r-2)$-space $\pi$ contained in $\mathcal{P}'_{r,e}$. We can use Equation~\eqref{eq:arbitrarypi} from 
Theorem~\ref{Th1Gen}, combined with the results from Lemma~\ref{Le:aid1} and~\ref{Le:aid2}. This yields the
following inequality. 

\begin{equation*}
\begin{split}
m(q^{e-1}+1)(m(q^{r+e-1}+1)-\mu(\pi))+q^{r+e-2}(m-\mu(\pi))(q^e{\color{red}+}1)-m(q^{r+e-2}+1)(m-\mu(\pi))(q^e+1)\\
\geq q^{e-1}\left( m(q^e+1)(m-\mu(\pi))+(1+\mu(\pi))(mq^e(q^{r-1}-1)+\mu(\pi)(q^e+1))\right)+\mu(\pi)q^{r+2e-1}\frac{q^{r-2}-1}{q-1}. 
\end{split}
\end{equation*}
This inequality then simplifies to
\begin{equation}
\begin{split}
m^2(q^{r+e-1}&-q^{r+e-2}-q^{2e-1}-q^e) \\
&+m\left( \mu(\pi)(q^{r+e-2}+2q^{2e-1}+q^e)+q^{r+e-2}+q^{2e-1}\right)\\
&-\mu(\pi)\left( q^{r+2e-2}+q^{r+e-2}+(1+\mu(\pi))(q^{2e-1}+q^{e-1})+q^{r+2e-1}\frac{q^{r-2}-1}{q-1}\right)\geq 0
\end{split}
\end{equation}
Dividing by $q^{e-1}$ gives the statement of the Lemma.
\end{proof}

\begin{St*}[Theorem \ref{th:MainImprov}]
Let $q > 2$ and $r \geq 3$.  Suppose that $\mathcal{O}$ is an $m$-ovoid in $\mathcal{P}_{r,e}'$, 
with (a) $r\geq 4$, or, (b) $e\in\{1,\frac32\}$ and $(r,q,e)\not= (3,3,1)$. Then it holds that
$$m\geq\frac{-r(1+\frac{2}{q^{r-e-1}}+\frac{1}{q^{r-2}})+\sqrt{r^2(1+\frac{1}{q^{r-e-1}})^2+4(q-2)(r-1)(q^{e+1}\frac{q^{r-2}-1}{q-1}+q^e+1)}}{2(q-1)}.$$
This bound asymptotically converges to $$m\geq \frac{-r+\sqrt{r^2+4(r-1)(q-2)q^{r+e-2}}}{2(q-1)}.$$ 
\end{St*}
\begin{proof}
Suppose that $\pi$ is an arbitrary $(r-2)$-dimensional subspace contained in $\mathcal{P}_{r,e}'$, 
such that $\mu(\pi^\perp\setminus \pi)\not=0$. Then we can apply Theorem \ref{th:MainImprovV1}, and
we denote Inequality~\eqref{EqNew} as 
$$\alpha m^2+\beta m +\delta \geq 0$$
where $\alpha$, $\beta$ and $\delta$ depend on $r,e,q$ and $\mu(\pi)$. Clearly, if $q > 2$ and $r \geq 3$, then 
$\alpha > 0$. Hence it follows that
$$m\geq \frac{-\beta+\sqrt{\beta^2-4\alpha\delta}}{2\alpha}.$$
We distinguish two cases. 

First assume that $m> r-1$. Then by Lemma~\ref{Le:aid3}, we can find an 
$(r-2)$-dimensional subspace $\pi$ contained in $\mathcal{P}_{r,e}'$ for which 
$\mu(\pi)=\min\{r-1, m\}=r-1$. Since 
$$\mu(\pi^\perp\setminus \pi)=(m-\mu(\pi))(q^{e}+1)\not=0.$$ 
we can apply Theorem \ref{th:MainImprovV1} on the subspace $\pi$.

To approximate the lower bound on $m$, we will use an upper bound for $\alpha$ in denominator, a lower bound for 
$\alpha$ in the numerator, an upper bound for $\beta$ in the numerator, and a lower bound for $\beta^2$ in the numerator. As such,
we will obtain an easier to handle expression. 

We find that 
\begin{itemize}
    \item $\alpha \leq q^{r-1}(q-1),$
    \item $r(q^{r-1}+q^{e})\leq \beta \leq r(q^{r-1}+2q^{e}+q),$ %since $2q^{e-1}\leq q^{2e-1}$ 
    for $q>2$ and using that $\mu(\pi) = r-1$.
\end{itemize}
Finally, if (a)  $r\geq 4$ or (b) $r\geq 3$ and $e\in\{1,\frac32\}$, we find that $q^{r-1}\geq q^{e}+q$. 
Thus $\alpha\geq q^{r}-2q^{r-1}$. Hence we find that
    \begin{equation}\label{eq:theorem5.5}
    \begin{split}
        \alpha\cdot (-\delta) & \geq (r-1)\left( q^{r+e}\frac{q^{r-2}-1}{q-1}+q^{r+e-1}+q^{r-1}\right)\alpha\\
        &\geq(r-1) q^{r-1}\left( q^{e+1}\frac{q^{r-2}-1}{q-1}+q^{e}+1\right)q^{r-1}(q-2).
    \end{split}
    \end{equation}
From this we can conclude that
\begin{equation}\label{eq:bound1}
    m\geq\frac{-r(1+\frac{2}{q^{r-e-1}}+\frac{1}{q^{r-2}})+\sqrt{r^2(1+\frac{1}{q^{r-e-1}})^2+4(q-2)(r-1)(q^{e+1}\frac{q^{r-2}-1}{q-1}+q^e+1)}}{2(q-1)}.
\end{equation}

Now assume that $m\leq r-1$. Under the condition that (a) $r\geq 4$, or, (b) $e\in\{1,\frac32\}$ and $(r,q,e)\not= (3,3,1)$, it is 
easy to check that $r-1 < b$, with $b$ the lower bound for $m$ of Theorem~\ref{ThSmallImprov}. Hence the 
assumption $m \leq r-1$ leads immediately to a contradiction. 

\end{proof}

We can adapt the approximation in the proof of the previous theorem slightly to handle the case $r=3$, $e=2$.

\begin{St}\label{th:mainQ7}
    Suppose that $\mathcal{O}$ is an $m$-ovoid in $\q^-(7,q)$, for $q>2$, then 
    $$ m\geq\frac{-3(3+\frac{1}{q})+\sqrt{36+8(q-\frac73)(q^{3}+q^2+1)}}{2(q-1)}.$$
\end{St}
\begin{proof}
    The proof is completely similar to the proof of Theorem \ref{th:MainImprov}. 
    To find a lower bound on $m$ using Inequality~\eqref{eq:theorem5.5}, we need some adjustments. 
    Note that $\frac43q^{2}\geq q^{2}+q$, which results in $\alpha\geq q^{3}-\frac73q^{2}$ and thus
  $$\alpha\cdot (-\delta) \geq 2 q^{2}\left( q^{3}+q^{2}+1\right)q^{2}(q-\frac{7}{3}).$$
Hence
\begin{equation}
    m\geq\frac{-3(3+\frac{1}{q})+\sqrt{9(1+1)^2+8(q-\frac73)(q^{3}+q^2+1)}}{2(q-1)}.
\end{equation}
The rest of the proof follows completely similarly. In particular, the case for $m\leq r-1$ yields again a contradiction with Theorem \ref{ThSmallImprov}.
\end{proof}

\begin{Opm}\label{Re:improvements}
Comparing Theorem \ref{th:MainImprov} with Theorem \ref{ThSmallImprov}, it gives an improvement for 
many cases for $(r,q,e)$. In Tables \ref{tableW}, \ref{tableQ} and Table \ref{tableH} one can find the 
comparison between both results for $(r,3,1), (r,3,2)$ and $(r,9,3/2)$, with $r\leq 7$. Clearly in 
almost all cases the new bound yields an improvement. Secondly, observing the asymptotic case in 
Theorem \ref{th:MainImprov}, it clearly yields an improvement from Theorem \ref{ThSmallImprov}. 
This is also shown in Tables \ref{tableW}, \ref{tableQ} and Table \ref{tableH} for $r=100$.    
Generally the improvement enlarges whenever $r$ and $q$ are bigger. Finally, 
the improvement that arises from Theorem \ref{th:mainQ7} can be seen in Table \ref{tableQ7}.
\end{Opm}

\begin{table}[h]
\begin{center}
\begin{tabular}{|c|c|c|}
\hline
$r$ & Bound from Theorem \ref{th:MainImprov} & Bound from Theorem \ref{ThSmallImprov} \\
\hline
%3 & / & $m\geq 2$ \\
%\hline
4 &$m\geq 5$ & $m\geq 4$ \\
\hline
5 & $m\geq 10$ & $m\geq 8$ \\
\hline
6 & $m\geq 20$ & $m\geq 13$ \\
\hline
7 & $m\geq 39$ & $m\geq 23$ \\
\hline
100 & $m\geq 2{.}53 \times 10^{24}$ & $m\geq 3{.}59 \times 10^{23}$\\
\hline
\end{tabular}
\caption{Bounds for $m$-ovoids of $W(2r-1,3)$.}\label{tableW}
\end{center}
\end{table}

\begin{table}[h]
\begin{center}
\begin{tabular}{|c|c|c|}
\hline
$r$ & Bound from Theorem \ref{th:MainImprov} & Bound from Theorem \ref{ThSmallImprov} \\
%\hline
%3 &  $m\geq2$ & $m\geq 4$ \\
\hline
4 & $m\geq 8$ & $m\geq 8$ \\
\hline
5 & $m\geq 18$ & $m\geq 13$ \\
\hline
6 & $m\geq 36$ & $m\geq 23$ \\
\hline
7 & $m\geq 69$ & $m\geq 40$ \\
\hline
100 & $m\geq 4{.}37 \times 10^{24}$ & $m\geq 6{.}22 \times 10^{23}$\\
\hline
\end{tabular}
\caption{Bounds for $m$-ovoids of $Q^-(2r+1,3)$.}\label{tableQ}
\end{center}
\end{table}

\begin{table}[h]
\begin{center}
\begin{tabular}{|c|c|c|}
\hline
$r$ & Bound from Theorem \ref{th:MainImprov} & Bound from Theorem \ref{ThSmallImprov} \\
\hline
3 & $m\geq 8$ & $m\geq 6$ \\
\hline
4 & $m\geq 29$ & $m\geq 18$ \\
\hline
5 & $m\geq 99$ & $m\geq 53$ \\
\hline
6 & $m\geq 330$ & $m\geq 158$ \\
\hline
7 & $m\geq 1085$ & $m\geq 474$ \\
\hline
100 & $m\geq 1{.}04\times 10^{48}$ & $m\geq 1{.}12\times 10^{47}$\\
\hline
\end{tabular}
\caption{Bounds for $m$-ovoids of $H(2r,3^2)$.}\label{tableH}
\end{center}
\end{table}

\begin{table}[h]
\begin{center}
\begin{tabular}{|c|c|c|}
\hline
$q$ & Bound from Theorem \ref{th:mainQ7} & Bound from Theorem \ref{ThSmallImprov} \\
\hline
3 & $m\geq 2$ & $m\geq 4$ \\
\hline
4 &$m\geq 4$ & $m\geq 5$ \\
\hline
5 &$m\geq 6$ & $m\geq 6$ \\
\hline
7 & $m\geq 10$ & $m\geq 8$ \\
\hline
8 & $m\geq 11$ & $m\geq 9$ \\
\hline
$3^5$ & $m\geq 345$ & $m\geq 244$\\
\hline
\end{tabular}
\caption{Bounds for $m$-ovoids of $Q^-(7,q)$.}\label{tableQ7}
\end{center}
\end{table}

\section*{Acknowledgments}
The authors thank the anonymous referees for their careful reading and valuable comments and 
Francesco Pavese for the fruitful discussions.

%\bibliographystyle{plain}
%\bibliography{ref}

\end{document}